\documentclass[11pt]{amsart}
\usepackage{amsthm}
\usepackage{amssymb}
\usepackage{amsmath}
\newtheorem{lemma}{Lemma}[section]
\newtheorem{theorem}[lemma]{Theorem}
\newtheorem*{theorem*}{Theorem}

\newtheorem{proposition}{Proposition}[section]

\newtheorem{corollary}{Corollary}[section]
\newtheorem*{corollary*}{Corollary}

\numberwithin{equation}{section}

\gdef\myletter{}

\let\savetheequation\theequation
\def\theequation{\savetheequation\myletter}

\newcommand{\CC}{{\mathbb C}}

\newcommand{\RR}{{\mathbb R}}

\def \bar{\overline}

\begin{document}

\vskip 3mm

\title[Transfinite diameter]{\bf Transfinite diameter notions in $\CC^N$ and integrals of Vandermonde determinants}  

\author{T. Bloom* and N. Levenberg}{\thanks{*Supported in part by NSERC grant}}

\address{University of Toronto, Toronto, Ontario M5S 2E4 Canada}  
\email{bloom@math.toronto.edu}

\address{Indiana University, Bloomington, IN 47405 USA}

\email{nlevenbe@indiana.edu}

\date{December 
15, 2007}

\begin{abstract}
We provide a general framework and indicate relations between the notions of transfinite diameter, homogeneous transfinite diameter, and weighted transfinite diameter for sets in $\CC^N$. An ingredient is a formula of Rumely \cite{rumely} which relates the Robin function and the transfinite diameter of a compact set. We also prove limiting formulas for integrals of generalized Vandermonde  determinants with varying weights for a general class of compact sets and measures in $\CC^N$. Our results extend to certain weights and measures defined on cones in $\RR^N$.
\end{abstract}

\maketitle

\section{\bf Introduction.}
	\label{sec:intro}
	
	Given a compact set $E$ in the complex plane $\CC$, the {\it transfinite diameter} of $E$ is the number
	$$d(E):=\lim_{n\to\infty} \max_{\zeta_1,...,\zeta_n\in E}|VDM(\zeta_1,...,\zeta_n)|^{1/{n \choose 2}}:=\max_{\zeta_1,...,\zeta_n\in E}\prod_{i<j}|\zeta_i-\zeta_j|^{1/{n \choose 2}}.$$
	It is well-known that this quantity is equivalent to the {\it Chebyshev constant} of $E$:
	$$T(E):=\lim_{n\to \infty} [\inf\{||p_n||_E: p_n(z)=z^n +\sum_{j=1}^{n-1}c_jz^j\}]^{1/n}$$
	(here, $||p_n||_E:=\sup_{z\in E}|p(z)|$) and also to $e^{-\rho(E)}$ where
	$$\rho(E):= \lim_{|z|\to \infty} [g_E(z) -\log {|z|}]$$
	is the {\it Robin constant} of $E$. The function $g_E$ is the {\it Green function} of logarithmic growth associated to $E$. Moreover, if $w$ is an admissible weight function on $E$, weighted versions of the above quantities can be defined. We refer the reader to the book of Saff-Totik \cite{safftotik} for the definitions and relationships. 
	
	For $E\subset \CC^N$ with $N>1$, multivariate notions of transfinite diameter, Chebyshev constant and Robin-type constants have been introduced and studied by several people. For an introduction to weighted versions of some of these quantitites, see Appendix B by Bloom in \cite{safftotik}. In the first part of this paper (section 2), we discuss a general framework for the various types of transfinite diameters in the spirit of Zaharjuta \cite{zah}. In particular, we relate (Theorem 2.7) two weighted transfinite diameters, $d^w(E)$ and $\delta^w(E)$, of a compact set $E\subset \CC^N$ using a remarkable result of Rumely \cite{rumely} which itself relates the (unweighted) transfinite diameter $d(E)$ with a Robin-like integral formula. Very recently Berman-Boucksom \cite{BB} have established a generalization of Rumely's formula which includes a weighted version of his result.  
	
	In the second part of the paper (section 3) we generalize to $\CC^N$ some results on strong asymptotics of Christoffel functions proved in \cite{bloomlev} in one variable. For $E$ a compact subset of $\CC$, $w$ an admissible weight function on $E$, and $\mu$ a positive Borel measure on $E$ such that the triple $(E,w,\mu)$ satisfies a weighted Bernstein-Markov inequality (see (\ref{wtdbernmark})), we take, for each $n=1,2,...$, a set of orthonormal polynomials $q_1^{(n)},...,q_{n}^{(n)}$ with respect to the varying measures $w(z)^{2n}d\mu(z)$ where deg$q_j^{(n)}=j-1$ and form the sequence of Christoffel functions $K_n(z):=\sum_{j=1}^{n} |q_j^{(n)}(z)|^2$. In \cite{bloomlev} we showed that 
	\begin{equation} \label{thm22}\frac{1}{n}K_n(z)w(z)^{2n}d\mu(z)  \to d\mu_{eq}^w(z)\end{equation} 
	weak-* where $\mu_{eq}^w$ is the potential-theoretic weighted equilibrium measure. The key ingredients to proving (\ref{thm22}) are, firstly, the verification that 
	\begin{equation} \label{znasymp}\lim_{n\to \infty} Z_n^{1/n^2} = \delta^w(E)\end{equation} where
\begin{equation} \label{wtdzn1} Z_n=Z_n(E,w,\mu):= \end{equation}
$$\int_{E^n} |VDM(\lambda_1,...,\lambda_n)|^2w(\lambda_1)^{2n} \cdots w(\lambda_n)^{2n}d\mu(\lambda_1) \cdots d\mu(\lambda_n);$$
and, secondly, a ``large deviation'' result in the spirit of Johansson \cite{joh}. We generalize these two results to $\CC^N, \ N>1$ (Theorems 3.1 and 3.2). The methods are similar to the corresponding one variable methods and were announced in \cite{bloomlev}, Remark 3.1. In particular, $\delta^w(E)$ is interpreted as the transfinite diameter of a circled set in one higher dimension. We also discuss the case where $E=\Gamma$ is an unbounded cone 
in $\RR^N$ for special weights and measures. Some of our results were proved independently by Berman-Boucksom \cite{BB}. 

We end the paper with a short section which includes questions related to these topics. We are grateful to Robert Berman for making reference \cite{BB} available to us. The second author would also like to thank Sione Ma'u and Laura DeMarco for helpful conversations.
	
	\section{\bf Transfinite diameter notions in $\CC^N$.}
	\label{sec:gf}

	We begin by considering a function $Y$ from the set of multiindices $\alpha \in {\bf N}^N$ to the nonnegative real numbers satisfying:
\begin{equation} \label{subadd} Y(\alpha + \beta) \leq Y(\alpha)\cdot Y(\beta) \ \hbox{for all} \ \alpha, \ \beta \in {\bf N}^N. \end{equation}
We call a function $Y$ satisfying (\ref{subadd}) {\it submultiplicative};  we have three main examples below. Let $e_1(z),...,e_j(z),...$ be a listing of the monomials
$\{e_i(z)=z^{\alpha(i)}=z_1^{\alpha_1}\cdots z_N^{\alpha_N}\}$ in
$\CC^N$ indexed using a lexicographic ordering on the multiindices $\alpha=\alpha(i)=(\alpha_1,...,\alpha_N)\in {\bf N}^N$, but with deg$e_i=|\alpha(i)|$ nondecreasing. We write $|\alpha|:=\sum_{j=1}^N\alpha_j$. 

We define the following integers:
\begin{enumerate}
\item $m_d^{(N)}=m_d:=$ the number of monomials $e_i(z)$ of degree at most $d$ in $N$ variables;
\item $h_d^{(N)}=h_d:=$ the number of monomials $e_i(z)$ of degree exactly $d$ in $N$ variables;
\item $l_d^{(N)}=l_d:=$ the sum of the degrees of the $m_d$ monomials $e_i(z)$ of degree at most $d$ in $N$ variables.
\end{enumerate}

\noindent We have the following relations:

\begin{equation} \label{precrucial} m_d^{(N)} ={N+d \choose d}; \ h_d^{(N)}=m_d^{(N)}-m_{d-1}^{(N)}= {N-1+d \choose d} \end{equation}
and
\begin{equation} \label{crucial} h_d^{(N+1)} =  {N+d \choose d}= m_d^{(N)}; \ l_d^{(N)} = N{N+d \choose N+1}=(\frac{N}{N+1})\cdot dm_d^{(N)}. \end{equation}
The elementary fact that the dimension of the space of homogeneous polynomials of degree $d$ in $N+1$ variables equals the dimension of the space of polynomials of degree at most $d$ in $N$ variables will be  crucial in sections 4 and 5. Finally, we let
$$r_d^{(N)}=r_d:=dh_d^{(N)}=d(m_d^{(N)}-m_{d-1}^{(N)})$$
which is the sum of the degrees of the $h_d$ monomials $e_i(z)$ of degree exactly $d$ in $N$ variables. We observe that
\begin{equation} \label{postcrucial} l_d^{(N)}=\sum_{k=1}^dr_k^{(N)}=\sum_{k=1}^N kh_k^{(N)}.\end{equation}

Let $K\subset \CC^N$ be compact. Here are our three natural constructions associated to $K$:

\begin{enumerate}
\item {\it Chebyshev constants}: Define the class of polynomials $$P_i=P(\alpha(i)):=\{e_i(z)+\sum_{j<i}c_je_j(z)\};  $$
and the Chebyshev constants
$$Y_1(\alpha):= \inf \{||p||_K:p\in P_i\}.$$ 
We write $t_{\alpha,K}:=t_{\alpha(i),K}$ for a Chebyshev polynomial; i.e., $t_{\alpha,K} \in P(\alpha(i))$ and $||t_{\alpha,K}||_K=Y_1(\alpha)$.  
\item {\it Homogeneous Chebyshev constants}: Define the class of homogeneous polynomials $$P_i^{(H)}=P^{(H)}(\alpha(i)):=\{e_i(z)+\sum_{j<i, \ {\rm deg}(e_j)={\rm deg}(e_i)}c_je_j(z)\};$$
and the homogeneous Chebyshev constants
$$Y_2(\alpha):= \inf \{||p||_K:p\in P^{(H)}_i\}.$$
We write $t^{(H)}_{\alpha,K}:=t^{(H)}_{\alpha(i),K}$ for a homogeneous Chebyshev polynomial; i.e., $t^{(H)}_{\alpha,K} \in P^{(H)}(\alpha(i))$ and $||t^{(H)}_{\alpha,K}||_K=Y_2(\alpha)$.
\item {\it Weighted Chebyshev constants}: Let $w$ be an admissible weight function on $K$ (see below) and let 
$$Y_3(\alpha):= \inf \{||w^{|\alpha (i)|}p||_K:=\sup_{z\in K}\{|w(z)^{|\alpha (i)|}p(z)|:p\in P_i\}$$ 
be the weighted Chebyshev constants. Note we use the polynomial classes $P_i$ as in (1). We write $t^w_{\alpha,K}$ for a weighted Chebyshev polynomial; i.e., $t^w_{\alpha,K}$ is of the form $w^{\alpha(i)}p$ with $p\in P(\alpha(i))$ and 
$||t^w_{\alpha,K}||_K= Y_3(\alpha)$.
\end{enumerate}

Let
$\Sigma$ denote the standard
$(N-1)-$simplex  in $\RR^N$; i.e.,
$$\Sigma = \{\theta =(\theta_1,...,\theta_N)\in \RR^N: \sum_{j=1}^N\theta_j=1, \ \theta_j\geq 0, \ j=1,...,N\},$$
and let
$$\Sigma^0 :=\{\theta \in \Sigma:  \ \theta_j > 0, \ j=1,...,N\}.$$
Given a submultiplicative function $Y(\alpha)$, define, as with the above examples, a new function  
\begin{equation} \label{tau}\tau( \alpha):= Y(\alpha)^{1/|\alpha|}.\end{equation} 
An examination of lemmas 1, 2, 3, 5, and 6 in \cite {zah} shows that (\ref{subadd}) is the only property of the numbers $Y( \alpha)$ needed to establish those lemmas. That is,
we have the following results for $Y:{\bf N}^N \to \RR^+$ satisfying (\ref{subadd}) and the associated function $\tau( \alpha)$ in (\ref{tau}).

\begin{lemma} For  all $\theta \in \Sigma^0$, the limit
$$T(Y,\theta):= \lim_{\alpha/|\alpha|\to \theta} Y(\alpha)^{1/|\alpha|}=\lim_{\alpha/|\alpha|\to \theta}\tau(\alpha)$$
exists. \end{lemma}

\begin{lemma} The function $\theta \to T(Y,\theta)$ is log-convex on $\Sigma^0$ (and hence continuous). \end{lemma}

\begin{lemma}  Given $b\in \partial \Sigma$,
$$\liminf_{\theta \to b, \ \theta \in \Sigma^0}
T(Y,\theta) =\liminf_{i\to \infty, \ \alpha(i)/|\alpha(i)|\to b}\tau(\alpha(i)).$$
\end{lemma}

\begin{lemma}  Let $\theta(k):= \alpha(k)/|\alpha(k)|$ for $k=1,2,...$ and let $Q$ be a compact subset of $\Sigma^0$. Then
$$\limsup_{|\alpha|\to \infty} \{\log \tau(\alpha(k)) - \log T (Y(\theta(k))): |\alpha(k)|=\alpha, \ \theta(k)\in Q\}=0.$$
\end{lemma}

\begin{lemma} \label{lemma2.5} Define 
$$\tau(Y):= \exp \bigl[\frac{1}{\hbox{meas}(\Sigma)}\int_{\Sigma} \log T(Y,\theta)d\theta\bigr]$$
Then
$$\lim_{d\to \infty} \frac{1}{h_d}\sum_{|\alpha|=d}\log \tau(\alpha)= \log \tau (Y);$$
i.e., using (\ref{tau}),
$$\lim_{d\to \infty}\bigl[\prod_{|\alpha|= d}Y(\alpha)\bigr]^{1/dh_d} = \ \tau (Y).$$
 \end{lemma}

One can incorporate all of the $Y(\alpha)'$s for $|\alpha|\leq d$; this is the content of the next result. 

\begin{theorem}\label{geommean} We have 
$$\lim_{d\to \infty}\bigl[\prod_{|\alpha|\leq d}Y(\alpha)\bigr]^{1/l_d} \ \hbox{exists and equals} \ \tau (Y).$$
\end{theorem}

\begin{proof} Define the geometric means 
$$\tau^0_d := \bigl(\prod_{|\alpha| =d}\tau(\alpha)\bigr)^{1/h_d}, \ d=1,2,...$$
The sequence
$$\log \tau^0_1,\log \tau^0_1,... (r_1 \ \hbox{times}), ... , \log \tau^0_d, \log \tau^0_d,... (r_d \ \hbox{times}), ...$$
converges to $\log \tau(Y)$ by the previous lemma; hence the arithmetic mean of the first $l_d = \sum_{k=1}^d r_k$ terms (see (\ref{postcrucial})) converges to $\log \tau(Y)$ as well. Exponentiating this arithmetic mean gives
\begin{equation} \label{mainest}\bigl(\prod_{k=1}^d (\tau^0_k)^{r_k} \bigr)^{1/l_d}= \bigl(\prod_{k=1}^d\prod_{|\alpha|=k}\tau(\alpha)^k\bigr)^{1/l_d}=\bigl(\prod_{|\alpha|\leq d}Y(\alpha)\bigr)^{1/l_d}\end{equation}
and the result follows.
\end{proof}

Returning to our examples (1)-(3), example (1) was the original setting of Zaharjuta \cite{zah} which he utilized to prove the existence of the limit in the definition of the {\it transfinite diameter} of a compact set $K\subset \CC^N$. For 
$\zeta_1,...,\zeta_n\in \CC^N$, let
\begin{equation} \label{vdm}VDM(\zeta_1,...,\zeta_n)=\det [e_i(\zeta_j)]_{i,j=1,...,n}  \end{equation}
$$= \det
\left[
\begin{array}{ccccc}
 e_1(\zeta_1) &e_1(\zeta_2) &\ldots  &e_1(\zeta_n)\\
  \vdots  & \vdots & \ddots  & \vdots \\
e_n(\zeta_1) &e_n(\zeta_2) &\ldots  &e_n(\zeta_n)
\end{array}
\right]$$
and for a compact subset $K\subset \CC^N$ let
$$V_n =V_n(K):=\max_{\zeta_1,...,\zeta_n\in K}|VDM(\zeta_1,...,\zeta_n)|.$$
Then 
\begin{equation} \label{tdlim}d(K)=\lim_{d\to \infty}V_{m_d}^{1/l_d} \end{equation} is the {\it transfinite diameter} of $K$; Zaharjuta \cite{zah} showed that the limit exists by showing that one has
\begin{equation} \label{zahthm}d(K)=\exp\bigl[\frac{1}{\hbox{meas}(\Sigma)}\int_{\Sigma^0}\log {\tau(K,\theta)}d\theta \bigr]  \end{equation} 
where $\tau(K,\theta)=T(Y_1,\theta)$ from (1); i.e., the right-hand-side of (\ref{zahthm}) is $\tau(Y_1)$. This follows from Theorem \ref{geommean} for $Y=Y_1$ and the estimate
$$\bigl(\prod_{k=1}^d (\tau^0_k)^{r_k} \bigr)^{1/l_d} \leq V_{m_d}^{1/l_d}   \leq (m_d!)^{1/l_d} \bigl(\prod_{k=1}^d (\tau^0_k)^{r_k} \bigr)^{1/l_d}$$
in \cite{zah} (compare (\ref{mainest})).

For a compact {\it circled} set $K\subset \CC^N$; i.e., $z\in K$ if and only if $e^{i\phi}z\in K, \ \phi \in [0,2\pi]$, one need only
consider homogeneous polynomials in the definition of the directional Chebyshev constants $\tau(K,\theta)$. In other words, in the notation of (1) and (2), $Y_1(\alpha)=Y_2(\alpha)$ for all $\alpha$ so that 
$$T(Y_1,\theta)=T(Y_2,\theta) \ \hbox{for circled sets} \  K.$$
This is because for such a set, if we write a polynomial $p$ of degree $d$ as $p=\sum_{j=0}^dH_j$ where $H_j$ is a homogeneous
polynomial of degree $j$, then, from the Cauchy integral formula, $||H_j||_K\leq ||p||_K, \ j=0,...,d$. Moreover, a slight modification of Zaharjuta's arguments prove the existence of the limit of appropriate roots of maximal {\it homogeneous} Vandermonde determinants; i.e., the homogeneous transfinite diameter $d^{(H)}(K)$ of a compact set (cf., \cite{jed}). From the above remarks, it follows that 
\begin{equation} \label{circled} \hbox{for circled sets} \  K, \  d(K)=d^{(H)}(K). \end{equation}
Since we will be using the homogeneous transfinite diameter, we amplify the discussion. We relabel the standard basis monomials $\{e_i^{(H,d)}(z)=z^{\alpha(i)}=z_1^{\alpha_1}\cdots z_N^{\alpha_N}\}$ where $|\alpha(i)|=d, \ i=1,...,h_d$, we define the $d-$homogeneous Vandermonde determinant 
\begin{equation} \label{vdmh}VDMH_d((\zeta_1,...,\zeta_{h_d}):=\det \bigl[e_i^{(H,d)}(\zeta_j)\bigr]_{i,j=1,...,h_d}.\end{equation}
Then 
\begin{equation} \label{htdlim}d^{(H)}(K)=\lim_{d\to \infty}\bigl[\max_{\zeta_1,...,\zeta_{h_d}\in K}|VDMH_d(\zeta_1,...,\zeta_{h_d})|\bigr]^{1/dh_d} \end{equation} is the homogeneous transfinite diameter of $K$; the limit exists and equals 
$$\exp\bigl[\frac{1}{\hbox{meas}(\Sigma)}\int_{\Sigma^0}\log {T(Y_2,\theta)}d\theta \bigr]$$
where $T(Y_2,\theta)$ comes from (2).

Finally, related to example (3), there are similar properties for the weighted version of directional Chebyshev constants and transfinite diameter. To define weighted notions, let $K\subset {\bf
C}^N$ be closed and let $w$ be an {\it admissible} weight function on $K$; i.e., $w$ is a nonnegative, usc function with
$\{z\in K:w(z)>0\}$. Let $Q:= -\log w$ and define the weighted
pluricomplex Green function $V^*_{K,Q}(z):=\limsup_{\zeta \to z}V_{K,Q}(\zeta)$ where
$$V_{K,Q}(z):=\sup \{u(z):u\in L(\CC^N), \ u\leq Q \ \hbox{on} \ K\}. $$
Here, $L(\CC^N)$ is the set of all plurisubharmonic functions $u$ on $\CC^N$ with the property that $u(z) - \log |z| = 0(1), \ |z| \to \infty$. If $K$ is closed but not necessarily bounded, we require that $w$ satisfies the growth property
\begin{equation} \label{grprop} |z|w(z)\to 0 \ \hbox{as} \ |z|\to \infty, \ z\in K, 
\end{equation}
so that $V_{K,Q}$ is well-defined and equals $V_{K\cap B_R,Q}$ for $R>0$ sufficiently large where $B_R=\{z:|z|\leq R\}$ (Definition 2.1 and Lemma 2.2 of Appendix B in \cite{safftotik}). The unweighted case is when $w\equiv 1$ ($Q\equiv 0$); we then write $V_K$ for the pluricomplex Green function. The set $K$ is called {\it regular} if $V_K=V_K^*$; i.e., $V_K$ is continuous; and $K$ is {\it locally regular} if for each $z\in K$, the sets $K\cap \overline{B(z,r)}$ are regular for $r>0$ where $B(z,r)$ denotes the ball of radius $r$ centered at $z$. We define the {\it weighted transfinite diameter}
$$d^w(K):=\exp\bigl[\frac{1}{\hbox{meas}(\Sigma)}\int_{\Sigma^0}\log {\tau^w(K,\theta)}d\theta \bigr] $$
as in \cite{bloomlev2} where $\tau^w(K,\theta)= T(Y_3,\theta)$ from (3); i.e., the right-hand-side of this equation is the quantity $\tau(Y_3)$.

We remark for future use that if $\{K_j\}$ is a decreasing sequence of locally regular compacta with $K_j \downarrow K$, and if $w_j$ is a continuous admissible weight function on $K_j$ with $w_j \downarrow w$ on $K$ where $w$ is an admissible weight function on $K$, then the argument in Proposition 7.5 of  \cite{bloomlev2} shows that $\lim_{j\to \infty} \tau^{w_j}(K_j,\theta) = \tau^w(K,\theta)$ for all $\theta \in \Sigma^0$ (we mention that there is a misprint in the statement of this proposition in \cite{bloomlev2}) and hence
\begin{equation}\label{decrlimit1} \lim_{j\to \infty} d^{w_j}(K_j) = d^w(K). \end{equation}
In particular, (\ref{decrlimit1}) holds in the unweighted case ($w\equiv 1$) for any decreasing sequence $\{K_j\}$ of compacta with $K_j \downarrow K$; i.e., 
\begin{equation}\label{decrlimit2} \lim_{j\to \infty} d(K_j) = d(K) \end{equation}
(cf., \cite{bloomlev2} equation (1.13)).

Another natural definition of a weighted transfinite diameter uses weighted Vandermonde determinants. Let $K\subset \CC^N$ be compact and let
$w$ be an admissible weight function on
$K$.  Given $\zeta_1,...,\zeta_n\in K$, let
$$W(\zeta_1,...,\zeta_n):=VDM(\zeta_1,...,\zeta_n)w(\zeta_1)^{|\alpha (n)|}\cdots w(\zeta_n)^{|\alpha (n)|}$$
$$= \det
\left[
\begin{array}{ccccc}
 e_1(\zeta_1) &e_1(\zeta_2) &\ldots  &e_1(\zeta_n)\\
  \vdots  & \vdots & \ddots  & \vdots \\
e_n(\zeta_1) &e_n(\zeta_2) &\ldots  &e_n(\zeta_n)
\end{array}
\right]\cdot w(\zeta_1)^{|\alpha (n)|}\cdots w(\zeta_n)^{|\alpha (n)|}$$
be a {\it weighted Vandermonde determinant}. Let
\begin{equation} \label{wn} W_n:=\max_{\zeta_1,...,\zeta_n\in K}|W(\zeta_1,...,\zeta_n)|
\end{equation}
and define an {\it $n-$th weighted Fekete set for $K$ and $w$} to be a set of $n$ points $\zeta_1,...,\zeta_ n\in K$ with the property that
$$|W(\zeta_1,...,\zeta_n)|=\sup_{\xi_1,...,\xi_n\in K}|W(\xi_1,...,\xi_n)|.$$
Also, define
\begin{equation} \label{deltaw} \delta^w(K):=\limsup_{d\to \infty}W_{m_d}^{1/ l_d}. \end{equation}

We will show in Proposition \ref{wtdlimit} that $\lim_{d\to \infty}W_{m_d}^{1/ l_d}$ (the weighted analogue of (\ref{tdlim})) exists. The question of the existence of this limit if $N>1$ was raised in \cite{bloomlev2}. Moreover, using a recent result of Rumely, we show how $\delta^w(K)$ is related to $d^w(K)$:
 \begin{equation} \label{conj1}\delta^w(K)=[\exp {(-\int_{K} Q(dd^cV^*_{K,Q})^{N})}]^{1/N}\cdot d^w(K)\end{equation}
where $(dd^cV^*_{K,Q})^{N}$ is the complex Monge-Ampere operator applied to $V^*_{K,Q}$. We refer the reader to \cite{klimek} or Appendix B of \cite{safftotik} for more on the complex Monge-Ampere operator.  
	
	We begin by proving the existence of the limit in the definition of $\delta^w(E)$ in (\ref{deltaw}) for a set $E\subset \CC^N$ and an admissible weight $w$ on $E$ (see also \cite{BB}).
		
	\begin{proposition} \label{wtdlimit} Let $E\subset \CC^N$ be a compact set with an admissible weight function $w$. The limit 
	$$\delta^w(E):=\lim_{d\to \infty} \bigl[\max_{\lambda^{(i)}\in E}|VDM(\lambda^{(1)},...,\lambda^{(m_d^{(N)})})|\cdot w(\lambda^{(1)})^d\cdots w(\lambda^{(m_d^{(N)})})^d\bigr] ^{1/l_d^{(N)}}$$
	exists. 
		\end{proposition}
		
		\begin{proof} Following \cite{bloom}, we define the circled set 
	$$F=F(E,w):=\{(t,z)=(t,t \lambda )\in \CC^{N+1}: \lambda \in E, \ |t|=w(\lambda)\}.$$
	We first relate weighted Vandermonde determinants for $E$ with homogeneous Vandermonde determinants for $F$. To this end, for each positive integer $d$, choose 
	$$m_d^{(N)}=  {N+d \choose d}$$
	(recall (\ref{precrucial})) points $\{(t_i,z^{(i)})\}_{i=1,...,m_d^{(N)}}=\{(t_i,t_i\lambda^{(i)})\}_{i=1,...,m_d^{(N)}}$ in $F$ and form the $d-$homogeneous Vandermonde determinant
	$$VDMH_d((t_1,z^{(1)}),...,(t_{m_d^{(N)}},z^{(m_d^{(N)})})).$$ 
	We extend the lexicographical order of the monomials in $\CC^N$ to $\CC^{N+1}$ by letting $t$ precede any of $z_1,...,z_N$. Writing the standard basis monomials of degree $d$ in 
	$\CC^{N+1}$ as 
	$$\{t^{d-j} e_{k}^{(H,d)}(z): j=0,...,d; \ k=1,...,h_j\};$$ 
	i.e., for each power $d-j$ of $t$, we multiply by the standard basis monomials of degree $j$ in $\CC^{N}$, and dropping the superscript $(N)$ in $m_d^{(N)}$, we have the $d-$homogeneous Vandermonde matrix
$$\left[\begin{array}{ccccc}
 t_1^d &t_2^d&\ldots  &t_{m_d}^d\\
 t_1^{d-1}e_{2}(z^{(1)})  &t_2^{d-1}e_{2}(z^{(2)})&\ldots  &t_{m_d}^{d-1}e_{2}(z^{(m_d)})\\
  \vdots  & \vdots & \ddots  & \vdots \\
e_{m_d}(z^{(1)}) &e_{m_d}(z^{(2)} ) &\ldots  &e_{m_d}(z^{(m_d)})
\end{array}\right]$$
$$=\left[\begin{array}{ccccc}
 t_1^d &t_2^d&\ldots  &t_{m_d}^d\\
 t_1^{d-1}z_1^{(1)}  &t_2^{d-1}z_1^{(2)}&\ldots  &t_{m_d}^{d-1}z_1^{(m_d)}\\
  \vdots  & \vdots & \ddots  & \vdots \\
(z^{(1)}_N)^d &(z^{(2)}_N)^d  &\ldots  &(z^{(m_d)}_N)^d
\end{array}\right].$$
	Factoring $t_i^d$ out of the $i-$th column, we obtain
	$$VDMH_d((t_1,z^{(1)}),...,(t_{m_d},z^{(m_d)}))=t_1^d\cdots t_{m_d}^d\cdot VDM(\lambda^{(1)},...,\lambda^{(m_d)});$$
thus, writing $|A|:=|\det A|$ for a square matrix $A$,  
\begin{align} \label{matrix}
\left|\begin{array}{ccccc}
 t_1^d &t_2^d&\ldots  &t_{m_d}^d\\
 t_1^{d-1}z_1^{(1)}  &t_2^{d-1}z_1^{(2)}&\ldots  &t_{m_d}^{d-1}z_1^{(m_d)}\\
  \vdots  & \vdots & \ddots  & \vdots \\
(z^{(1)}_N)^d &(z^{(2)}_N)^d  &\ldots  &(z^{(m_d)}_N)^d
\end{array}\right|\end{align}
$$=|t_1|^d\cdots |t_{m_d}|^d \left|\begin{array}{ccccc}
 1 &1 &\ldots  &1\\
 \lambda^{(1)}_1 &\lambda^{(2)}_1 &\ldots  &\lambda^{(m_d)}_1\\
  \vdots  & \vdots & \ddots  & \vdots \\
(\lambda^{(1)}_N)^d &(\lambda^{(2)}_N)^d &\ldots  &(\lambda^{(m_d)}_N)^d
\end{array}\right|,$$
where $\lambda^{(j)}_k= z^{(j)}_k/t_j$ provided $t_j\not = 0$.  By definition of $F$, since $(t_i,z^{(i)})=(t_i,t_i\lambda^{(i)})\in F$, we have $ |t_i|=w(\lambda^{(i)})$ so that from (\ref{matrix})
	$$VDMH_d((t_1,z^{(1)}),...,(t_{m_d},z^{(m_d)}))$$
	$$=VDM(\lambda^{(1)},...,\lambda^{(m_d)})\cdot w(\lambda^{(1)})^d\cdots w(\lambda^{(m_d)})^d.$$
	Thus
	$$\max_{(t_i,z^{(i)})\in F}|VDMH_d((t_1,z^{(1)}),...,(t_{m_d},z^{(m_d)}))|=$$
	$$\max_{\lambda^{(i)}\in E}|VDM(\lambda^{(1)},...,\lambda^{(m_d)})|\cdot w(\lambda^{(1)})^d\cdots w(\lambda^{(m_d)})^d.$$
	Note that the maximum will occur when all $t_j=w(\lambda^{(j)})>0$. As mentioned in section \ref{sec:mainres} the limit
	$$\lim_{d\to \infty} \bigl[\max_{(t_i,z^{(i)})\in F}|VDMH_d((t_1,z^{(1)}),...,(t_{m_d},z^{(m_d)}))|\bigr]^{1/dh_d^{(N+1)}}=:d^{(H)}(F)$$
	exists \cite{jed}; thus the limit 
	$$\lim_{d\to \infty} \bigl[\max_{\lambda^{(i)}\in E}|VDM(\lambda^{(1)},...,\lambda^{(m_d)})|\cdot w(\lambda^{(1)})^d\cdots w(\lambda^{(m_d)})^d\bigr] ^{1/l_d^{(N)}}:=\delta^w(E)$$
	exists. \end{proof}
	
	\begin{corollary}\label{deltaw} For $E\subset \CC^N$ a nonpluripolar compact set with an admissible weight function $w$ and $$F=F(E,w):=\{(t,z)=(t,t \lambda )\in \CC^{N+1}: \lambda \in E, \ |t|=w(\lambda)\},$$
		\begin{equation} \label{homvswtd} \delta^w(E)=d^{(H)}(F)^{\frac{N+1}{N}}=d(F)^{\frac{N+1}{N}}. \end{equation}
		\end{corollary}
		\begin{proof} The first equality follows from the proof of Proposition \ref{wtdlimit} using the relation 
	$$l_d^{(N)}=(\frac{N}{N+1})\cdot dh_d^{(N+1)}$$
	(see (\ref{crucial})). The second equality is (\ref{circled}). 
	\end{proof}
	
	We next relate $\delta^w(E)$ and $d^w(E)$ but we first recall the remarkable formula of Rumely \cite{rumely}. For a plurisubharmonic function $u$ in $L(\CC^N)$ we can define the
{\it Robin function} associated to $u$:
$$\rho_u(z):=\limsup_{|\lambda|\to \infty} \left[u(\lambda z)-
\log(|\lambda|)\right].$$
This function is plurisubharmonic (cf., \cite{bloomold}, Proposition 2.1) and logarithmically homogeneous: 
$$\rho_u(tz) = \rho_u(z) + \log |t| \ \hbox{for} \ t\in \CC.$$
For $u=V_{E,Q}^*$ ($V_E^*$) we write $\rho_u= \rho_{E,Q}$ ($\rho_E$). Rumely's formula relates $\rho_E$ and $d(E)$:
\begin{align} \label{rum}  -\log d(E) = \frac{1}{N}\bigl[\int_{\CC^{N-1}}\rho_E(1,t_2,...,t_N)(dd^c\rho_E(1,t_2,...,t_N))^{N-1}\end{align}
$$+ \int_{\CC^{N-2}}\rho_E(0,1,t_3,...,t_N)(dd^c\rho_E(0,1,t_3,...,t_N))^{N-2}$$
$$+\cdots + \int_{\CC}\rho_E(0,..,0,1,t_N)(dd^c\rho_E(0,..,0,1,t_N) + \rho_E(0,..,0,1)\bigr].$$
Here we make the convention that $dd^c = \frac{1}{2\pi}(2i \partial \bar \partial)$ so that in any dimension $d=1,2,...$, 
$$\int_{\CC^d} (dd^cu)^d =1$$
for any $u\in L^+(\CC^d)$; i.e., for any plurisubharmonic function $u$ in $\CC^d$ which satisfies 
$$C_1 +\log (1+|z|) \leq u(z) \leq C_2 + \log (1+|z|)$$ for some $C_1,C_2$. 

We begin by rewriting (\ref{rum}) for {\it regular circled} sets $E$ using an observation of Sione Ma'u. Note that for such sets, $V_E^*=\rho_E^+:=\max(\rho_E,0)$. If we intersect $E$ with a hyperplane ${\mathcal H}$ through the origin, e.g., by rotating coordinates, we take ${\mathcal H}=\{z=(z_1,...,z_N)\in \CC^N: z_1 =0\}$, then $E\cap {\mathcal H}$ is a regular, compact, circled set in $\CC^{N-1}$ (which we identify with ${\mathcal H}$). Moreover, we have
$$\rho_{{\mathcal H}\cap E} (z_2,...,z_N)=\rho_E(0,z_2,...,z_N)$$
since each side is logarithmically homogeneous and vanishes for $(z_2,...,z_N)\in \partial ({\mathcal H}\cap E)$. Thus the terms 
$$\int_{\CC^{N-2}}\rho_E(0,1,t_3,...,t_N)(dd^c\rho_E(0,1,t_3,...,t_N))^{N-2}$$
$$+\cdots + \int_{\CC}\rho_E(0,..,0,1,t_N)(dd^c\rho_E(0,..,0,1,t_N) + \rho_E(0,..,0,1)$$
in (\ref{rum}) are seen to equal 
$$(N-1)d^{\CC^{N-1}}({\mathcal H}\cap E)$$ (where we temporarily write $d^{\CC^{N-1}}$ to denote the transfinite diameter in $\CC^{N-1}$ for emphasis) by applying (\ref{rum}) in $\CC^{N-1}$ to the set ${\mathcal H}\cap E$. Hence we have
\begin{align}  \label{newrum}  -\log d(E) = \frac{1}{N}\int_{\CC^{N-1}}\rho_E(1,t_2,...,t_N)(dd^c\rho_E(1,t_2,...,t_N))^{N-1} \end{align}
$$+\bigl(\frac{N-1}{N}\bigr)[-\log d^{\CC^{N-1}}({\mathcal H}\cap E)].$$

	\begin{theorem} For $E\subset \CC^N$ a nonpluripolar compact set with an admissible weight function $w$,
	 \begin{equation} \label{conj1}\delta^w(E)=[\exp {(-\int_{E} Q(dd^cV^*_{E,Q})^{N})}]^{1/N}\cdot d^w(E).\end{equation}
	\end{theorem}
\begin{proof} We first assume that $E$ is locally regular and $Q$ is continuous. It is known in this case that $V_{E,Q}=V_{E,Q}^*$ (cf., \cite{siciak}, Proposition 2.16). As before, we define the circled set 
	$$F=F(E,w):=\{(t,z)=(t,t \lambda )\in \CC^{N+1}: \lambda \in E, \ |t|=w(\lambda)\}.$$
	We claim this is a regular set; i.e., $V_F$ is continuous. First of all, $V_F^*(t,z)=\max[\rho_F(t,z),0]$ (cf., Proposition 2.2 of \cite{bloom}) so that it suffices to verify that $\rho_F(t,z)$ is continuous. From Theorem 2.1 and Corollary 2.1 of \cite{bloom}, 
	\begin{equation} \label{relation}V_{E,Q}(\lambda)=  \rho_F(1,\lambda) \ \hbox{on} \ \CC^N\end{equation}  
which implies, by the logarithmic homogeneity of $\rho_F$, that $\rho_F(t,z)$ is continuous on $\CC^{N+1}\setminus \{t=0\}$. Corollary 2.1 and equation (2.8) in \cite{bloom} give that
\begin{equation} \label{relation2}\rho_F(0,\lambda)= \rho_{E,Q}(\lambda) \ \hbox{for} \ \lambda \in \CC^N \end{equation}
and $\rho_{E,Q}$ is continuous by Theorem 2.5 of \cite{bloomlev2}. Moreover, the limit exists in the definition of $\rho_{E,Q}$:
$$\rho_{E,Q}(\lambda):=\limsup_{|t|\to \infty} [V_{E,Q}(t\lambda)-\log |t|]=\lim_{|t|\to \infty} [V_{E,Q}(t\lambda)-\log |t|];$$
and the limit is uniform in $\lambda$ (cf., Corollary 4.4 of \cite{blm}) which implies, from (\ref{relation}) and (\ref{relation2}), that $\lim_{t\to 0} \rho_F(t,\lambda)=\rho_F(0,\lambda)$ so that $\rho_F(t,z)$ is continuous. 
In particular, 
$$V_{E,Q}(\lambda)=Q(\lambda) = \rho_F(1,\lambda) \ \hbox{on the support of} \  (dd^cV_{E,Q})^{N}$$ so that
\begin{equation} \label{qandv} \int_{E} Q(\lambda)(dd^cV_{E,Q}(\lambda))^{N}=\int_{\CC^{N}} \rho_F(1,\lambda)(dd^c\rho_F(1,\lambda))^{N}.\end{equation}
 On the other hand, $E^w_{\rho}:=\{\lambda\in \CC^{N}: \rho_{E,Q}(\lambda)\leq 0\}$ is a circled set, and, according to eqn. (3.14) in \cite{bloomlev2}, 
$d^w(E)=d(E^w_{\rho})$. But 
$$\rho_{E,Q}(\lambda)=\limsup_{|t|\to \infty} [V_{E,Q}(t\lambda)-\log |t|]$$
$$=\limsup_{|t|\to \infty} [\rho_F(1,t\lambda)-\log |t|]=\limsup_{|t|\to \infty} \rho_F(1/t,\lambda)=\rho_F(0,\lambda).$$
Thus
$$E^w_{\rho}=\{\lambda\in \CC^{N}: \rho_F(0,\lambda)\leq 0\}=F\cap {\mathcal H}$$
where ${\mathcal H}=\{(t,z)\in \CC^{N+1}: t=0\}$ and hence 
\begin{equation} \label{step1} d^w(E)=d(E^w_{\rho})=d(F\cap {\mathcal H}).\end{equation}
From (\ref{newrum}) applied to $F\subset \CC^{N+1}$,
\begin{equation} \label{step2} -\log d(F)=  \frac{1}{N+1}\int_{\CC^{N}}\rho_F(1,\lambda)(dd^c\rho_F(1,\lambda))^{N}\end{equation}
$$+ (\frac{N}{N+1})[-\log d(F\cap {\mathcal H})].$$
Finally, from (\ref{homvswtd}),
\begin{equation} \label{step3}  \delta^w(E)=d(F)^{\frac{N+1}{N}};\end{equation}
putting together (\ref{qandv}), (\ref{step1}),  (\ref{step2}) and (\ref{step3}) gives the result if $E$ is locally regular and $Q$ is continuous. 

The general case follows from approximation. Take a sequence of locally regular compacta $\{E_j\}$ decreasing to $E$ and a sequence of weight functions $\{w_j\}$ with $w_j$ continuous and admissible on $E_j$ and $w_j \downarrow w$ on $E$ (cf., Lemma 2.3 of \cite{bloom}).  From (\ref{decrlimit1}) we have
\begin{equation} \label{deltalim}
\lim_{j\to \infty} d^{w_j}(E_j) =d^w(E). 
\end{equation}
Also, by Corollary \ref{deltaw} we have 
\begin{equation} \label{cor41}\delta^{w_j}(E_j)=d(F_j)^{\frac{N+1}{N}} \end{equation} where
$$F_j=F_j(E_j,w_j)=\{(t(1,\lambda):\lambda\in E_j, \ |t|=w_j(\lambda)\}.$$
Since $E_{j+1}\subset E_j$ and $w_{j+1}\leq w_j$, the sets 
$$\tilde F_j=\tilde F_j(E_j,w_j)=\{(t(1,\lambda):\lambda\in E_j, \ |t|\leq w_j(\lambda)\}$$
satisfy $\tilde F_{j+1} \subset \tilde F_j$ and hence
$$d(\tilde F_{j+1}) = d(F_{j+1})\leq d(\tilde F_{j})=d(F_{j}).$$
Since $F_j\downarrow F$, we conclude from (\ref{decrlimit2}) and (\ref{cor41}) that
\begin{equation} \label{dlimitok}
\lim_{j\to \infty} \delta^{w_j}(E_j) =\delta^w(E). 
\end{equation}
Applying (\ref{conj1}) to $E_j, w_j, Q_j$ and using (\ref{deltalim}) and (\ref{dlimitok}), we conclude that 
$$\int_{E_j} Q_j(dd^cV_{E_j,Q_j})^{N}\to \int_{E} Q(dd^cV^*_{E,Q})^{N},$$ completing the proof of (\ref{conj1}).
\end{proof}
\smallskip

\vskip 5mm

\section{\bf Integrals of Vandermonde determinants.}
	\label{sec:mainres}
	In this section, we first state and prove the analogue of an ``unweighted'' generalization to $\CC^N$ of Theorem 2.1 of \cite{bloomlev} as it has a self-contained proof. We first recall some terminology. Given a compact set $E\subset \CC^N$ and a measure $\nu$ on $E$, we say that $(E,\nu)$ satisfies the Bernstein-Markov inequality for holomorphic polynomials in $\CC^N$ if, given $\epsilon >0$, there exists a constant $ M= M(\epsilon)$ such that for all such polynomials $Q_n$ of degree at most $n$
	\begin{equation} \label{bernmark} ||Q_n||_E\leq  M(1+\epsilon)^n||Q_n||_{L^2(\nu)}. \end{equation}

\begin{theorem} \label{zntheoremsub} Let $(E,\mu)$ satisfy a Bernstein-Markov inequality. Then
$$\lim_{d\to \infty} Z_d^{1/2l_d^{(N)}} = d(E)$$ where
\begin{equation} \label{znunwtd} Z_d=Z_d(E,\mu):= \end{equation}
$$\int_{E^{m_d^{(N)}}} |VDM(\lambda^{(1)},...,\lambda^{m_d^{(N)}})|^2 d\mu(\lambda^{(1)}) \cdots d\mu(\lambda^{m_d^{(N)}}).$$
\end{theorem}

\begin{proof} Since $VDM(\zeta_1,...,\zeta_n)=\det [e_i(\zeta_j)]_{i,j=1,...,n}$ for any positive integer $n$, if we apply the Gram-Schmidt procedure to the monomials $e_1,...,e_{m_d^{(N)}}$ to obtain orthogonal polynomials $q_1,...,q_{m_d^{(N)}}$ with respect to $\mu$ where $q_j \in P_j$ has minimal $L^2(\mu)-$norm among all such polynomials, we get, upon using elementary row operations on $VDM(\zeta_1,...,\zeta_{m_d^{(N)}})$ and expanding the determinant (cf.,  \cite{deift} Chapter 5 or section 2 of \cite{bloomlev}) 
\begin{equation} \label{znl2} \int_{E^{m_d^{(N)}}}|VDM(\zeta_1,...,\zeta_{m_d^{(N)}})|^2d\mu(\zeta_1) \cdots d\mu(\zeta_{m_d^{(N)}}) = m_d^{(N)}! \prod_{j=1}^{m_d^{(N)}} ||q_j||_{L^2(\mu)}^2.\end{equation}
Let $t_{\alpha,E}\in P(\alpha)$ be a Chebyshev polynomial; i.e., $||t_{\alpha,E}||_E = Y_1(\alpha)$. Then Theorem \ref{geommean}  shows that
$$\lim_{d\to \infty}\bigl(\prod_{|\alpha|\leq d}||t_{\alpha,E}||_E \bigr)^{1/l_d}=\tau(Y_1)$$
since 
$$\lim_{d\to \infty} (m_d^{(N)}!)^{1/l_d^{(N)}}=1.$$
Zaharjuta's theorem (\ref{zahthm}) shows that $\tau(Y_1)=d(E)$ so we need show that
\begin{equation} \label{minimality}\lim_{d\to \infty}\bigl(\prod_{|\alpha|\leq d}||t_{\alpha,E}||_E \bigr)^{1/l_d}=\lim_{d\to \infty}\bigl(\prod_{|\alpha|\leq d}||q_{\alpha}||_{L^2(\mu)} \bigr)^{1/l_d}.\end{equation}
This follows from the Bernstein-Markov property. First note that
$$||q_{\alpha}||_{L^2(\mu)}\leq ||t_{\alpha,E}||_{L^2(\mu)}\leq \mu(E)\cdot ||t_{\alpha,E}||_E$$
from the $L^2(\mu)-$norm minimality of $q_{\alpha}$; then, given $\epsilon >0$, the Bernstein-Markov property and the sup-norm minimality of $t_{\alpha,E}$ give
$$||t_{\alpha,E}||_E \leq ||q_{\alpha}||_E\leq M(1+\epsilon)^{|\alpha|}||q_{\alpha}||_{L^2(\mu)}$$
for some $M=M(\epsilon)>0$. Taking products of these inequalities over $|\alpha|\leq d$; taking $l_d-$th roots; and letting $\epsilon \to 0$ gives the result. This reasoning is adapted from the proof of Theorem 3.3 in \cite{BBCL}.
\end{proof}

		A weighted polynomial on $E$ is a function of the form $w(z)^np_n(z)$ where $p_n$ is a holomorphic polynomial of degree at most $n$. Let $\mu$ be a measure with support in $E$ such that $(E,w,\mu)$ satisfies a Bernstein-Markov inequality for weighted polynomials (referred to as a {\it weighted B-M inequality} in \cite{bloom}): given $\epsilon >0$, there exists a constant $M=M(\epsilon)$ such that for all weighted polynomials $w^np_n$
	\begin{equation} \label{wtdbernmark}||w^np_n||_E\leq M(1+\epsilon)^n||w^np_n||_{L^2(\mu)}.\end{equation}
		 
		 Generalizing Theorem \ref{zntheoremsub}, we have the following result.
	 
\begin{theorem} \label{zntheorem} Let $(E,w,\mu)$ satisfy a Bernstein-Markov inequality (\ref{wtdbernmark}) for weighted polynomials. Then
$$\lim_{d\to \infty} Z_d^{1/2l_d^{(N)}} = \delta^w(E)$$ where
\begin{equation} \label{wtdzn} Z_d=Z_d(E,w,\mu):= \int_{E^{m_d^{(N)}}} |VDM(\lambda^{(1)},...,\lambda^{(m_d^{(N)})})|^2\times \end{equation}
$$w(\lambda^{(1)})^{2d}\cdots w(\lambda^{(m_d^{(N)})})^{2d}d\mu(\lambda^{(1)}) \cdots d\mu(\lambda^{(m_d^{(N)})}).$$
\end{theorem}
		
	The proof of Theorem \ref{zntheorem} follows along the lines of section 3 of \cite{bloomlev}. Let $E\subset \CC^N$ be a nonpluripolar compact set with an admissible weight function $w$ and let $\mu$ be a measure with support in $E$ such that $(E,w,\mu)$ satisfies a Bernstein-Markov inequality for weighted polynomials. The integrand 
	$$|VDM(\lambda^{(1)},...,\lambda^{(m_d^{(N)})})|^2\cdot w(\lambda^{(1)})^{2d}\cdots w(\lambda^{(m_d^{(N)})})^{2d} $$ in the definition of $Z_d$ in (\ref{wtdzn}) has a maximal value on $E^{m_d^{(N)}}$ whose $1/2l_d^{(N)}$ root tends to $\delta^w(E)$. To show that the integrals themselves have the same property, we begin by constructing the circled set $F\subset \CC^{N+1}$ defined as in section 4:
	$$F=F(E,w):=\{(t,z)=(t,t \lambda )\in \CC^{N+1}: \lambda \in E, \ |t|=w(\lambda)\}.$$
	We construct a measure $\nu$ on $F$ associated to $\mu$ such that $(F,\nu)$ satisfies the Bernstein-Markov property for holomorphic polynomials in $\CC^{N+1}$; i.e., (\ref{bernmark}) holds. Define 
	$$\nu:= m_{\lambda} \otimes \mu, \ \lambda \in E$$
	where $m_{\lambda}$ is normalized Lebesgue measure on the circle $|t|=w(\lambda)$ in the complex $t-$plane given by  
	$$C_{\lambda}:=\{(t,t\lambda)\in \CC^{N+1}: t\in \CC\}.$$
	That is, if $\phi$ is continuous on $F$, 
	$$\int_F \phi(t,z) d\nu (t,z) = \int_E\bigl[\int_{C_{\lambda}}\phi(t,t\lambda)dm_{\lambda}(t)\bigr] d\mu(\lambda).$$
	Equivalently, if $\pi:\CC^{N+1} \to \CC^N$ via $\pi(t,z)=z/t:=\lambda$, then $\pi_*(\nu) =\mu$. The fact that $(F,\nu)$ satisfies the Bernstein-Markov property follows from Theorem 3.1 of \cite{bloom}. Moreover, if $p_1(t,z)$ and $p_2(t,z)$ are two homogeneous polynomials in $\CC^{N+1}$ of degree $d$, say, and we write $$p_j(t,z)=p_j(t,t \lambda )=t^dp_j(1,\lambda)=:t^d G_j(\lambda), \ j=1,2$$
	for univariate $G_j$, then it is straightforward to see that 
	\begin{equation}
\label{orthog}	\int_F p_1(t,z)\overline {p_2(t,z)}d\nu(t,z) = \int_E G_1(\lambda)\overline {G_2(\lambda)}w(\lambda)^{2d}d\mu(\lambda) \end{equation}
	(cf., \cite{bloom}, Lemma 3.1 and its proof). Note that if $$p(t,z)=t^iz^{\alpha}=t^iz_1^{\alpha_1} \cdots z_N^{\alpha_N}$$ with $|\alpha|=\alpha_1+\cdots + \alpha_N=d-i$, then 
	$$p(t,z)=t^d(z/t)^{\alpha}=t^d G(\lambda)=t^d\cdot \lambda_1^{\alpha_1} \cdots \lambda_N^{\alpha_N}$$
	where $G(\lambda)= \lambda^{\alpha}=\lambda_1^{\alpha_1} \cdots \lambda_N^{\alpha_N}$.

	\begin{proposition} \label{prop1} Let 
	$$\tilde Z_d:= \int_{F^{m_d^{(N)}}} |VDMH_d((t_1,z^{(1)}),...,(t_{m_d^{(N)}},z^{(m_d^{(N)})}))|^2$$
	$$ d\nu(t_1,z^{(1)}) \cdots d\nu(t_{m_d^{(N)}},z^{(m_d^{(N)})}).$$
	Then $\tilde Z_d= Z_d$ where $m_d^{(N)}=  {N+d \choose d}$ and (recall (\ref{wtdzn}))
	$$Z_d=\int_{E^{m_d^{(N)}}} |VDM(\lambda^{(1)},...,\lambda^{(m_d^{(N)})})|^2 \times $$
	$$ w(\lambda^{(1)})^{2d}\cdots w(\lambda^{(m_d^{(N)})})^{2d}d\mu(\lambda^{(1)}) \cdots d\mu(\lambda^{(m_d^{(N)})}).$$
 \end{proposition}
	
	\begin{proof} Recall from section 2 that the $d-$homogeneous Vandermonde determinant $VDMH_d((t_1,z^{(1)}),...,(t_{m_d^{(N)}},z^{(m_d^{(N)})}))$ equals 
$$\det \left[\begin{array}{ccccc}
 t_1^d &t_2^d&\ldots  &t_{m_d^{(N)}}^d\\
  \vdots  & \vdots & \ddots  & \vdots \\
e_{m_d^{(N)}}(z^{(1)}) &e_{m_d^{(N)}}(z^{(2)} ) &\ldots  &e_{m_d^{(N)}}(z^{(m_d^{(N)})})
\end{array}\right].$$
Expanding this determinant in $\tilde Z_d$ gives
	$$\tilde Z_d =\sum_{I,S}\sigma(I)\cdot \sigma(S) \bigl[\int_F t_1^{d-deg(e_{i_1})}e_{i_1}(z^{(1)})\bar t_1^{d-deg(e_{s_1})} \overline { e_{s_1}(z^{(1)})}d\nu(t_1,z^{(1)})\cdots $$
$$\cdots \int_F t_{m_d^{(N)}}^{d-deg(e_{i_{m_d^{(N)}}})}e_{i_{m_d^{(N)}}}(z^{(m_d^{(N)})})\bar t_{m_d^{(N)}}^{d-deg(e_{s_{m_d^{(N)}}})} \overline { e_{s_{m_d^{(N)}}}(z^{(m_d^{(N)})})}d\nu(t_{m_d^{(N)}},z^{(m_d^{(N)})})\bigr]$$
where $I=(i_1,...,i_{m_d^{(N)}})$ and $S=(s_1,...,s_{m_d^{(N)}})$ are permutations of $(1,...,m_d^{(N)})$ and $\sigma(I)$ is the sign of $I$ ($+1$ if $I$ is even; $-1$ if $I$ is odd). 
Expanding the ordinary Vandermonde determinant in $Z_d$ gives 
$$Z_d = \sum_{I,S}\sigma(I)\cdot \sigma(S) \bigl[\int_E e_{i_1}(\lambda^{(1)})\overline {e_{s_1}(\lambda^{(1)})}w(\lambda^{(1)})^{2d}d\mu(\lambda^{(1)})\cdots $$
$$\cdots \int_Ee_{i_{m_d^{(N)}}}(\lambda^{(m_d^{(N)})})\overline {e_{s_{m_d^{(N)}}}(\lambda^{(m_d^{(N)})})}w(\lambda^{(m_d^{(N)})})^{2d}d\mu(\lambda^{(m_d^{(N)})})\bigr].$$
Since $|t_j|=w(\lambda^{(j)})$, using (\ref{orthog}) completes the proof.
	\end{proof}
	
	We need to work in $\CC^{N+1}$ with the $\tilde Z_d$ integrals and verify the following.
	
	\begin{proposition} \label{prop2} We have 
	$$\lim_{n\to \infty} \tilde Z_d^{1/2dm_d^{(N)}}= d^{(H)}(F).$$	\end{proposition}
	
	\begin{proof} Fix $d$ and consider the $m_d^{(N)}$ monomials 
	$$t^d,  t^{d-1}z_1 , \cdots,  z_N^{d},$$ 
	utilized in 
	$\tilde VDMH_d((t_1,z^{(1)}),...,(t_{n},z^{(m_d^{(N)})}).$ Use Gram-Schmidt in $L^2(\nu)$ to obtain orthogonal homogeneous polynomials 
	$$q^{(H)}_1(t,z)=t^d, \ q^{(H)}_2(t,z) = t^{d-1}z_1 + ..., \cdots, q^{(H)}_{m_d^{(N)}}(t,z) = z_N^{d} + ....$$
	Then 
	$$\tilde VDMH_d((t_1,z^{(1)}),...,(t_{m_d^{(N)}},z^{(m_d^{(N)})})=\det \bigl[ q^{(H)}_i(t_j,z^{(j)})\bigr]_{i,j=1,...,m_d^{(N)}}.$$
	By orthogonality, as in (\ref{znl2}), we obtain
	$$\tilde Z_d =m_d^{(N)}! ||q^{(H)}_1||_{L^2(\nu)}^2\cdots ||q^{(H)}_{m_d^{(N)}}||_{L^2(\nu)}^2.$$
	Note that from (\ref{precrucial}) and (\ref{crucial}) $(m_d^{(N)}!)^{1/2dm_d^{(N)}}\to 1$ as $d\to \infty$. Now from Lemma 
	\ref{lemma2.5} we have 
	$$\lim_{d\to \infty}\bigl(\prod_{|\alpha|=d}||t^{(H)}_{\alpha,F}||_F\bigr)^{1/dm_d^{(N)}}=\tau(Y_2)=\tau(Y_1)=d^{(H)}(F).$$
	Thus we need to show that
	$$\lim_{d\to \infty}\bigl(\prod_{|\alpha|=d}||t^{(H)}_{\alpha,F}||_F\bigr)^{1/dm_d^{(N)}}=\lim_{d\to \infty}\bigl(\prod_{i=1}^{m_d^{(N)}}||q^{(H)}_i||_{L^2(\nu)}\bigr)^{1/dm_d^{(N)}}.$$
	This is analogous to (\ref{minimality}) in the proof of Theorem \ref{zntheoremsub} and it follows in the same manner from the Bernstein-Markov property for $(F,\nu)$ and the minimality properties of $t^{(H)}_{\alpha,F}$ and $q^{(H)}_i$.
	\end{proof}
	
	Combining Propositions \ref{prop1} and \ref{prop2} with equation (\ref{homvswtd}) and the second equation in (\ref{crucial}) completes the proof of Theorem \ref{zntheorem}. \hfill $\Box$
	
	As a corollary, we get a ``large deviation'' result, which follows easily from Theorem \ref{zntheorem}. Define a probability measure ${\mathcal P}_d$ on $E^{m_d^{(N)}}$ via, for a Borel set $A\subset E^{m_d^{(N)}}$,
	$${\mathcal P}_d(A):=\frac{1}{Z_d}\int_{A} |VDM(z_1,...,z_{m_d^{(N)}})|^2w(z_1)^{2d} \cdots w(z_{m_d^{(N)}})^{2d}d\mu(z_1) \cdots d\mu(z_{m_d^{(N)}}).$$
		
\begin{proposition} Given $\eta >0$, define
$$A_{d,\eta}:=\{(z_1,...,z_{m_d^{(N)}})\in E^{m_d^{(N)}}: $$
$$|VDM(z_1,...,z_{m_d^{(N)}})|^2w(z_1)^{2d} \cdots w(z_n)^{2d} \geq (\delta^w(E) -\eta)^{2l_d}\}.$$
Then there exists $d^*=d^*(\eta)$ such that for all $d>d^*$, 
$${\mathcal P}_d(E^{m_d^{(N)}}\setminus A_{d,\eta})\leq (1-\frac{\eta}{2\delta^w(E)})^{2l_d}.$$
	\end{proposition}	
	
	\begin{proof} From Theorem \ref{zntheorem}, given $\epsilon >0$, 
	$$Z_d \geq [\delta^w(E) -\epsilon]^{2l_d}$$
	for $d\geq d(\epsilon)$. Thus
	$${\mathcal P}_d(E^{m_d^{(N)}}\setminus A_{d,\eta})=$$
	$$\frac{1}{Z_d}\int_{E^{m_d^{(N)}}\setminus A_{d,\eta}} |VDM(z_1,...,z_{m_d^{(N)}})|^2w(z_1)^{2d} \cdots w(z_n)^{2d}d\mu(z_1) \cdots d\mu(z_{m_d^{(N)}})$$
	$$\leq \frac{  [\delta^w(E) -\eta]^{2l_d}} { [\delta^w(E) -\epsilon]^{2l_d}}$$
	if $d\geq d(\epsilon)$. Choosing $\epsilon < \eta/2$ and $d^*=d(\epsilon)$ gives the result.
	\end{proof}

Finally, we state a a version of (\ref{znasymp}) for $\Gamma$ an unbounded cone in $\RR^N$ with $\Gamma =\overline {\hbox{int}\Gamma}$. Precisely, our set-up is the following. Let $R(x)=R(x_1,...,x_N)$ be a polynomial in $N$ (real) variables $x=(x_1,...,x_N)$ and let 
\begin{equation} \label{definemu}d\mu(x):=|R(x)|dx = |R(x_1,...,x_N)|dx_1 \cdots dx_N. \end{equation} Next, let $w(x)=\exp(-Q(x))$ where $Q(x)$ satisfies the inequality
\begin{equation} \label{defineq} Q(x) \geq c |x|^{\gamma}\end{equation}
for all $x\in \Gamma$ for some $c,\gamma>0$. 

\begin{theorem} \label{znrn} Let $S_w:=$supp$(dd^cV_{\Gamma,Q})^N$ where $Q$ is defined as in (\ref{defineq}). With $\mu$ defined in (\ref{definemu}), 
$$\lim_{d\to \infty} Z_d^{1/2l_d^{(N)}} = \delta^w(S_w)$$ where
\begin{equation} \label{znrndef} Z_d=Z_d(\Gamma,w,\mu):=\int_{\Gamma^{m_d^{(N)}}} |VDM(\lambda^{(1)},...,\lambda^{(m_d^{(N)})})|^2\times \end{equation}
$$w(\lambda^{(1)})^{2d}\cdots w(\lambda^{(m_d^{(N)})})^{2d}d\mu(\lambda^{(1)}) \cdots d\mu(\lambda^{(m_d^{(N)})}).$$
\end{theorem}

\smallskip

\noindent {\bf Remark}. The integrals considered in Theorem \ref{znrn} may be considered as multivariate versions (i.e., with a multivariable Vandermonde determinant in the integrand rather than a one-variable Vandermonde determinant) of integrals of the form
$$\int_{\RR^{d}} |VDM(\lambda_1,...,\lambda_d)|^2e^{-dQ(\lambda_1)}\cdots e^{-dQ(\lambda_d)}d\lambda_1 \cdots d\lambda_d$$
considered in \cite{deift}, Chapter 6, arising in the joint probability distribution of eigenvalues of certain random matrix ensembles. They are also multivariate versions of Selberg integrals of Laguerre type (cf., \cite{mehta}, equation (17.6.5)) which, after rescaling by a factor of $d$, are of the form, for $\Gamma =[0,\infty)\subset \RR$ and $\alpha >0$, 
$$\int_{\Gamma^{d}} |VDM(\lambda_1,...,\lambda_d)|^2e^{-d\lambda_1}\cdots e^{-d\lambda_d}(\prod_{j=1}^d\lambda_j^{\alpha})d\lambda_1 \cdots d\lambda_d.$$

\begin{proof} We begin by observing that 
\begin{equation} \label{integrand}VDM(\lambda^{(1)},...,\lambda^{(m_d^{(N)})})^2\cdot w(\lambda^{(1)})^{2d}\cdots w(\lambda^{(m_d^{(N)})})^{2d}\end{equation}
	(the integrand in (\ref{znrndef}) with the absolute value removed from the VDM) becomes, if all but one of the $m_d^{(N)}-1$ variables are fixed, a weighted polynomial in the remaining variable. Since $w(x)$ is continuous, by Theorem 2.6 in Appendix B of \cite{safftotik}, a weighted polynomial attains its maximum on $S_w\subset \Gamma$. Hence the maximum value of (\ref{integrand}) on $\Gamma^{m_d^{(N)}}$ is attained on $(S_w)^{m_d^{(N)}}$. Since $S_w$ has compact support (cf., Lemma 2.2 of Appendix B of \cite{safftotik}), we can take $T>0$ sufficiently large with $S_w\subset \Gamma \cap B_T$ where $B_T:=\{x\in \RR^N: |x| \leq T\}$ and 
	$$\delta^w(S_w)=\delta^w(\Gamma \cap B_T).$$
	We need the following result.
	\begin{lemma} For all $T>0$ sufficiently large, there exists $M=M(T)>0$ with 
	$$||w^dp||_{L^2(\Gamma,\mu)}\leq M ||w^dp||_{L^2(\Gamma \cap B_T,\mu)}$$
	if $p=p(x)$ is a polynomial of degree $d$.
	\end{lemma}
	
	\begin{proof}  By Theorem 2.6 (ii) in Appendix B of \cite{safftotik}, we have
	$$|w(x)^dp(x)|\leq ||w^dp||_{S_w}e^{d(V_{\Gamma,Q}(x)-Q(x))}$$
	for all $x\in \Gamma$. Since $V_{\Gamma,Q}\in L(\CC^N)$ and $Q(x) \geq c |x|^{\gamma}$ for $x\in \Gamma$, there is a $c_0>0$ with 
	$$|w(x)^dp(x)|\leq ||w^dp||_{S_w}e^{-c_0d|x|^{\gamma}}$$
	for all $x\in \Gamma\cap B_T$ for $T$ sufficiently large. Hence
	$$||w^dp||_{L^2(\Gamma,\mu)}\leq ||w^dp||_{L^2(\Gamma \cap B_T,\mu)}
	+||w^dp||_{S_w}\int_{\Gamma \cap\{|x|\geq T\}}e^{-c_0d|x|^{\gamma}}|R(x)|dx.$$
	Now $(\Gamma \cap B_T,\mu)$ satisfies the Bernstein-Markov property (\cite{bloomlevold}, Theorem 2.1); thus by \cite{bloom} Theorem 3.2, the triple
	$(\Gamma\cap B_T,w,\mu)$ satisfies the weighted Bernstein-Markov property. Thus, given $\epsilon >0$, there is $M_1=M_1(\epsilon)>0$ with 
	$$||w^dp||_{S_w}=||w^dp||_{\Gamma \cap B_T}\leq M_1(1+\epsilon)^d ||w^dp||_{L^2(\Gamma \cap B_T,\mu)}.$$
	A simple estimate shows that 
	$$\int_{|x|\geq T}e^{-c_0d|x|^{\gamma}}|R(x)|dx \leq e^{-c'd}$$
	for some $c'>0$. The result now follows by choosing $\epsilon$ sufficiently small. 
	\end{proof}

	We now expand the integrands in the formulas for $Z_d(\Gamma):=Z_d(\Gamma,w,\mu)$ and in $Z_d(\Gamma \cap B_T):=Z_d(\Gamma \cap B_T,w|_{\Gamma \cap B_T},\mu|_{\Gamma \cap B_T})$ as a product of $L^2-$norms of orthogonal polynomials as in (\ref{znl2}), and then proceed as in the proof of Corollary 2.1 in section 5 of \cite{bloomlev} to conclude that
	$$\lim_{d\to \infty} Z_d(\Gamma)^{1/2l_d^{(N)}} =\lim_{d\to \infty} Z_d(\Gamma \cap B_T)^{1/2l_d^{(N)}} =\delta^w(\Gamma \cap B_T) =\delta^w(S_w).$$\end{proof}
	
\smallskip

\vskip 5mm

\section{\bf Final remarks and questions.}
	\label{sec:finrem}
	
	In this section we discuss some further results in the literature and pose some questions. Recall from section 2 that a $d-$th weighted Fekete set for $E\subset \CC^N$ and an admissible weight $w$ on $E$ is a set of $m_d$ points $\zeta_1^{(d)},...,\zeta_{m_d}^{(d)}\in E$ with the property that
$$|W(\zeta_1,...,\zeta_{m_d})|=\sup_{\xi_1,...,\xi_{m_d}\in E}|W(\xi_1,...,\xi_{m_d})|$$ where $W$ is defined in (\ref{wn}). In \cite{bloomlev2} the authors asked if the sequence of probability measures
$$\mu_d := {1\over m_d}\sum_{j=1}^{m_d} <\zeta_j^{(d)}>, \ d=1,2,...,$$
where $<z>$ denotes the point mass at $z$ and $\{\zeta_1^{(d)},...,\zeta_{m_d}^{(d)}\}$ is a 
$d-$th weighted Fekete set for $E$ and $w$, has a unique weak-* limit, and, if so, whether this limit is the Monge-Ampere measure, $\mu^w_{eq}:=(dd^cV_{E,Q}^*)^N$. From the proof of Proposition \ref{wtdlimit}, a $d-$th weighted Fekete set for $E$ and $w$ corresponds to a $d-th$ homogeneous Fekete set for the circled set 
$$F=F(E,w):=\{(t,z)=(t,t \lambda )\in \CC^{N+1}: \lambda \in E, \ |t|=w(\lambda)\};$$
i.e., a set of $m_d^{(N)}=h_d^{(N+1)}$ points in $F$ which maximize the corresponding homogeneous Vandermonde determinant (\ref{vdmh}) for $F$. From Theorem 2.2 of \cite{bloom}, to verify this conjecture for $E\subset \CC^N$ and an admissible weight $w$ it suffices to verify it for homogeneous Fekete points associated to circled sets in $\CC^{N+1}$.

\medskip	

	Suppose now that $\mu$ is a measure on $E$ such that $(E,w,\mu)$ satisfies a Bernstein-Markov inequality for weighted polynomials. Define the probability measures $$\mu_d(z):=\frac {1}{Z_d} R_1^{(d)}(z)w(z)^{2d}d\mu(z)$$ where $Z_d$ is defined in (\ref{wtdzn}) and
\begin{equation} \label{r1} R_1^{(d)}(z):=\int_{E^{m_d-1}} |VDM(\lambda^{(1)},...,\lambda^{(m_d-1)},z)|^2\cdot \end{equation}
$$w(\lambda^{(1)})^{2d} \cdots w(\lambda^{(m_d-1)})^{2d}d\mu(\lambda^{(1)}) \cdots d\mu(\lambda^{(m_d-1)}).$$
We observe that with the notation in (\ref{r1}) and (\ref{wtdzn})
\begin{equation}
\label{r1form} \frac {R_1^{(d)}(z)}{Z_d} = \frac{1}{m_d} \sum_{j=1}^{m_d} |q_j^{(d)}(z)|^2 \end{equation}
where $q^{(d)}_1,...,q^{(d)}_{m_d}$ are orthonormal polynomials with respect to the measure $w(z)^{2d}d\mu(z)$ forming a basis for the  polynomials of degree at most $d$. To verify (\ref{r1form}), we refer the reader to the argument in Remark 2.1 of \cite{bloomlev}. Forming the sequence of Christoffel functions $K_d(z):=\sum_{j=1}^{m_d} |q_j^{(d)}(z)|^2$, in \cite{bloomlev} Theorem 2.2 it was shown that if $N=1$ then  $\mu_d(z) \to \mu^w_{eq}(z)$ weak-*; i.e., 
	\begin{equation} \label{thmd22}\frac{1}{m_d}K_d(z)w(z)^{2d}d\mu(z) \to \mu_{eq}^w(z) \ \hbox{weak-* }.\end{equation} 
We conjecture that (\ref{thmd22}) should hold in $\CC^N$ for $N>1$. To this end, we remark that if $E=\bar D$ where $D$ is a smoothly bounded domain in $\RR^N$, it follows from the proof of Theorem 1.3 in \cite{bt} that $\mu_{eq}:=(dd^cV_{E}^*)^N=c(x)dx$ is absolutely continuous with respect to $\RR^N-$Lebesgue measure $dx$ on $D$; and if $\mu(x)=f(x)dx$ is also absolutely continuous, then a conjectured version of (\ref{thmd22}) in the unweighted case $w\equiv 1$ is
$$\frac{1}{m_d}K_d(x)f(x) \to c(x) \ \hbox{on} \ D.$$ 
Bos (\cite{bos1}, \cite{bos2}) has verified this result for centrally symmetric functions $f(x)$ on the unit ball in $\RR^N$ and Xu \cite{xu} proved this result for certain Jacobi-type functions $f(x)$ on the standard simplex in $\RR^N$. For further results on subsets of $\RR^1$, see references [10], [13], [15]-[17] in \cite{bloomlev}. Berman (\cite{Ber} and \cite{Ber2}) has shown that if $w=e^{-Q}$ is a smooth admissible weight function on $\CC^N$ (recall (\ref{grprop})), then $\mu_{eq}^w:=(dd^cV_{\CC^N,Q}^*)^N=c(z)dz$ is absolutely continuous with respect to $\CC^N-$Lebesgue measure $dz$ on the interior $I$ of the compact set $\{z\in \CC^N: V_{\CC^N,Q}(z) = Q(z)\}$ and 
$$\frac{1}{m_d}K_d(z)Q(z)^d \to c(z) \ \hbox{a.e. on} \ I.$$

\end{document}